\documentclass [12pt,a4paper,draft]{article}
\usepackage{latexsym,amsfonts,amssymb,amsmath,longtable,amsthm}
\usepackage{cite}
\newtheorem{lem}{Lemma}
\newtheorem{theorem}{Theorem}

\renewcommand{\ge}{\geqslant}

\newtheorem{Lemma}{Lemma}

\theoremstyle{definition}
\newtheorem{probl}[Lemma]{Problem}

\title{\vspace{-1cm} \hfill{\normalsize MSC2010: 20D06, 20E45, 20G41}{
\fontfamily{cmr} \fontseries{bx} \selectfont \\ \vspace{1cm} Strong reality of finite simple groups}
\thanks{The work is partially supported by  RFBR, projects 08-01-00322, 10-01-00391, and 10-01-90007,  ADTP ``Development of the Scientific
Potential of
Higher
School'' of the Russian Federal Agency for Education (Grant
2.1.1.419), Federal Target Grant "Scientific and
educational personnel of innovation Russia" for 2009-2013 (government
contract No. 02.740.11.0429). The first author gratefully acknowledges the support from Deligne 2004 Balzan prize in mathematics, and the
Lavrent'ev Young Scientists Competition (No 43 on 04.02.2010)}}
\date{}
\author{\bf  E.P.Vdovin,
A.A.Gal't}

\begin{document}

\sloppy

\maketitle
\pagenumbering{arabic}

\begin{abstract}
The classification of finite simple strongly real groups is complete. It is
easy to see that strong reality for
every nonabelian
finite simple group is equivalent to the fact that each element can be written
as a product of two involutions. We thus obtain a solution
to Problem 14.82 from the Kourovka notebook from  the classification
of finite simple strongly real groups.
\end{abstract}

\section*{Introduction}
In this article we solve Problem 14.82 from the Kourovka notebook [1].

\begin{probl}
{\rm [1, 14.82]} Find all finite simple groups whose every element 
is a product of two involutions.
\end{probl}

Since each involution of a nonabelian finite simple group lies in an elementary
abelian subgroup of order $4$; i.e., for every
involution $t$ there exist an involution $s\not=t$ commuting with $t$, Problem
14.82 is equivalent to that of the classification of
finite simple strongly real groups. Recall that an element $x$ of $G$ is called
{\em real} 
({\em strongly real}), if $x$ and 
$x^{-1}$ are conjugate in $G$ (respectively are conjugate by an involution in
$G$). A group $G$ is called {\em
real} ({\em strongly real}), if all elements of $G$ are real
(strongly real). Thus, if the order of $x$ is not
equal to  $1$ or $2$, then $x$ can be written as a product of two
involutions $s$ and $t$ if and only if 
$x$ is strongly real. Indeed, if $t$ is an involution inverting $x$ and $\vert
x\vert>2$ then $t$ and $tx$ are involutions and $x=t\cdot tx$.
Conversely, if there exists involutions $s$ and $t$ with $x=st$, then
$x^t=ts=x^{-1}$;
i.e., $x$ is strongly real.
The fact that in finite simple groups each element of order at most $2$ can be
written as a product of two involutions follows from the
Feit-Thompson Odd Order Theorem and the note  above. 

The problem of reality and strong reality of finite simple groups and the groups
in some sense close to simple was studied by many
authors, see [2--13]. In particular, the classification of finite simple real
groups is obtained in [2]. Thus it suffices to check what
finite simple real groups are strongly real in order to solve Problem 14.82. All
strongly real alternating and sporadic groups are found respectively in 
[3,\,4]. It is proven in  [5--7] that the symplectic groups
$\operatorname{PSp}_{2n}(q)$ are strongly real if and only if 
$q\not\equiv3\pmod4$. The strong reality of $\Omega_{4n}^{\varepsilon}(q)$ for
$q$ even is proven in [8]. The strong reality of
${\mathrm P}\Omega_{4n}^{-}(q)$ for $q$ odd is proven in [9]. Moreover, [10,
Theorem~8.5] implies  that if  $q$ is odd, then  ${\mathrm
P}\Omega_{4n}^{+}(q)$ and $\Omega_{2n+1}(q)$, together with ${\mathrm
P}\Omega_{4n}^{-}(q)$, are strongly real if $q\equiv1\pmod4$, while
$\Omega_9(q)$ and ${\mathrm P}\Omega_8^+(q)$ are also strongly real
if $q\equiv3\pmod4$. In the
present paper the following is proven.

\begin{theorem} {\em (Main Theorem)}
$G={}^3D_4(q)$ is strongly real.
\end{theorem}

\noindent The theorem together with 
[2--10]
imply the following theorems.

\begin{theorem}
Each finite simple real group is strongly real.
\end{theorem}

\begin{theorem}
Every element of a finite simple group $G$ can be written as a product of two
involutions if and only if $G$ is isomorphic to one of the
following groups:
\begin{itemize}
 \item [{\rm (1)}] $\operatorname{PSp}_{2n}(q)$ for $q\not\equiv3\pmod4$,
$n\ge1$;
\item[{\rm (2)}] $\Omega_{2n+1}(q)$ for $q\equiv1\pmod4$, $n\ge3$;
\item[{\rm (3)}] $\Omega_9(q)$ for  $q\equiv3\pmod4$;
\item[{\rm (4)}] ${\mathrm P}\Omega_{4n}^-(q)$ for $n\ge 2$;
\item[{\rm (5)}] ${\mathrm P}\Omega_{4n}^+(q)$ for  $q\not\equiv3\pmod4$,
$n\ge3$;
\item[{\rm (6)}] ${\mathrm P}\Omega_8^+(q)$;
\item[{\rm (7)}] ${}^3D_4(q)$;
\item[{\rm (8)}] $A_{10}$, $A_{14}$, $J_1$, $J_2$.
\end{itemize}
\end{theorem}
 
Theorem 3 gives a complete solution to Problem 14.82 from the Kourovka
notebook.

\section{Preliminary results}

Our notation for finite groups agrees with that of [14]. The notation and basic
facts for finite groups of Lie type and linear
algebraic groups can be found in [15]. A finite group $G$ is said to be a {\it
central product} of subgroups
$A$ and $B$ (which is denoted by $A\circ  B$) if $G=AB$ and the derived
subgroup $[A,B]$ is trivial. The order of a group $G$ and of
an element $g\in G$ we denote by  $| G|$ and $| g|$. If $X$ is a
subset of $G$ and 
$H$ is a subgroup of $G$ then the centralizer of $X$ in $G$ and the normalizer
of $H$ in $G$ are denoted by $C_G(X)$ and $N_G(H)$,
respectively. Given a subset $X$ of $G$ by $\langle X\rangle$  we denote the
subgroup
generated by $X$. A finite field of order $q$ we
denote
by ${\mathbb F}_q$, while $p$ always denotes its characteristic; i.e.,
$q=p^\alpha$
for some positive integer~$\alpha$.
By $e$ we denote the identity element of a group, while $1$ stands for  the unit
of a field.

Let $\overline{G}$ be a simple connected algebraic group over the algebraic
closure $\overline{{\mathbb F}}_p$ of a finite field 
${\mathbb F}_p$. A surjective endomorphism $\sigma$ of $\overline{G}$ is called
a
{\it  Steinberg endomorphism} (see [15,
Definition~1.15.1]), if the set of $\sigma$-stable points $\overline{G}_\sigma$
is finite.  $O^{p'}(\overline{G}_\sigma)$ is known
to be a finite group of Lie type, and each finite group of Lie type can be
obtained in this way (notice that given a finite group of Lie type a
corresponding algebraic group and a Steinberg map are not uniquely determined in
general). 
More detailed definitions and related results can be found in [15,
Sections~1.5,~2.2]. If $\overline{G}$ is simply connected then
$\overline{G}_\sigma=O^{p'}(\overline{G}_\sigma)$ by
[15, Theorem~2.2.6(f)].
Moreover, 
[16, Proposition~2.10] implies that the centralizer of every semisimple element
is a connected reductive subgroup of maximal rank
in~$\overline{G}$.

If $G$ is isomorphic to ${}^3D_4(q)$, then a corresponding algebraic group
$\overline{G}$ can be chosen simply connected. We always
assume that $\overline{G}$ is simply connected in this case, i.e., for every
$G={}^3D_4(q)$ a simply connected connected simple linear
algebraic group 
$\overline{G}=D_4(\overline{{\mathbb F}}_q)$, where $\overline{{\mathbb F}}_q$
is the
algebraic closure of
${\mathbb F}_q$, and a Steinberg endomorphism $\sigma$ are chosen so that 
$G=\overline{G}_\sigma$. In particular, the
centralizer of every semisimple element in $\overline{G}$ is connected. If
$\overline{T}$ is a 
$\sigma$-stable maximal torus of $\overline{G}$ then  $T=\overline{T}\cap G$ is
called a  {\it  maximal torus} of a finite group of Lie type
$G$. If $\overline{R}\leq \overline{S}$ are $\sigma$-stable subgroups of
$\overline{G}$, $R=\overline{R}\cap G$, and  $S=\overline{S}\cap
G$; then  $N_{\overline{S}}(\overline{R})\cap G$ is denoted by $N(S,R)$. Notice
that $N(S,R)\leq N_S(R)$, but the equality is not true in
general. For every $x\in G$ there exist unique elements $s,u\in G$ such that 
$x=su=us$, $s$ is semisimple, and  $u$ is unipotent. Furthermore, $s$ is the
$p'$-part of $x$, while  $u$ is the  $p$-part of 
$x$. This is called the {\it Jordan decomposition} of~$x$.

By [9, Lemma~10], all semisimple elements of ${}^3D_4(q)$ are strongly real.
Moreover, the conjugating involution found in the proof
of [9, Lemma~10] satisfies to the following property.

\begin{lem} For every maximal torus $T$ of ${}^3D_4(q)$ there exists an
involution
$x \in N(G,T)$ such that $t^x=t^{-1}$ for every ~${t \in T}$. In particular, for
every $t\in T$, both $xt$ and $tx$ are involutions
inverting every element of~$T$.
\end{lem}

All statements from the next lemma are immediate from the structure of
projective linear groups of degree~$2$.

\begin{lem}
The following hold:
  \begin{itemize}
 \item[{\em (1)}] $\operatorname{PSL}_2(q)$ is strongly real if and only if
$q\not\equiv3\pmod4$.
\item[{\em (2)}] $\operatorname{PGL}_2(q)$ is strongly real.
\item[{\em (3)}] If $q$ is odd, $u$ is a nonidentity unipotent element of 
$\operatorname{PGL}_2(q)$, and $t\in\operatorname{PGL}_2(q)$ is chosen so that
$u^t=u^k$ for some $k\in{\mathbb N}$; then $t$ lies in a Cartan
subgroup {\rm (}which is cyclic of order $q-1)$ of $\operatorname{PGL}_2(q)$
normalizing a unique maximal unipotent subgroup
of~$\operatorname{PGL}_2(q)$ containing $u$.
 \end{itemize}
\end{lem}

\section{Proof of the main theorem}

Let $G={}^3D_4(q)$ and $g\in G$. If $g$ is semisimple then by [9, Lemma~10] it
is strongly real. If
$g$ is unipotent and  $q$ is even then [12, Theorem~1] implies that $g$ is
strongly real. Assume that $g$ is unipotent, $q$ is odd and
$C_G(g)$ does not contain nonidentical semisimple elements; i.e., $C_G(g)$ is a
$p$-group. By [2, Lemma~5.9] there exists 
$x \in G$ such that $g^x=g^{-1}$.  Clearly we may assume that  $|x|=2^k$ for
some $k \in {\mathbb N}$. Then $x^2 \in C_G(g)$ and
$| x^2|$ is a power of~$2$. Therefore, $x^2$ is semisimple, whence $x^2=e$. If
$g$ is unipotent, $q$ is odd
and $C_G(g)$ contains a nonidentical semisimple element $s$; then we consider
$g_1=sg$. Decomposition $sg$ is the Jordan decomposition of
$g_1$. If we show the existence of an involution $x$ inverting $g_1$ then the
uniqueness of the Jordan decomposition implies 
$s^x=s^{-1}$ and $g^x=g^{-1}$. Thus we may assume that $g$ has a ``mixed
order''; i.e., in the Jordan decomposition
$g=su$ both $s$ and $u$ are nonidentical.

Assume that $C=C_G(s)$. Then $u \in C$. Moreover, 
$\overline{C}=C_{\overline{G}}(s)$ is a connected reductive subgroup of maximal
rank in
$\overline{G}$ and $C=\overline{C}_\sigma$. Clearly, every maximal torus $T$ of
$G$, which contains $s$, is included in $C$.
The structure of the centralizers of semisimple elements is given in [17,
Proposition~2.2]. Tables~2.2a and~2.2b from [17] are
the main technical instrument in the forthcoming arguments. If $q$ is even then
up to conjugation in $G$ there exist $8$
centralizers of order divisible by $p$ of nonidentical semisimple elements. If
$q$ is odd then there exist $9$ centralizers of this sort. We
consider each centralizer separately. Note that 
$\overline{C}=\overline{M}\circ\overline{S}$, where
$\overline{M}=[\overline{C},\overline{C}]$ is connected and semisimple, while 
$\overline{S}=Z(\overline{C})^0$ is a torus. Furthermore, $C$ possesses a
normal subgroup 
$M\circ S$, where $S=\overline{S}_\sigma\leq Z(C)$ and
$M=\overline{M}_\sigma=O^{p'}(C)$, and the structure of $M$ ($=M_{\sigma}$ in
the
notation from [17]) and $S$ ($=S_\sigma$ in the notation from [17]) is given
in [17, Tables~2.2a,~2.2b], where the structure or the order
of~${C/(M\circ S)}$ is also given. The indeces of elements below are chosen as
in [17, Tables~2.2a, ~2.2b]. Moreover, the subgroups and factor
groups of $C$ are isomorphic to classical groups in a natural way, and we
identify the subgroups and the factor groups of $C$ with the
corresponding classical groups.

Let $s$ be such that its centralizer is conjugate to the centralizer of $s_2$
(hence $s$ is an involution and this case can occur only if 
$q$ is odd). Then $M\simeq{\operatorname{SL}_2(q^3)}\circ
\operatorname{SL}_2(q)$, $|Z(M)|=2$ and $S=\{e\}$. Moreover, $|C:M|=2$ and, by
[18, Theorem~2], 
$C/\operatorname{SL}_2(q)\simeq \operatorname{PGL}_2(q^3)$ and
$C/\operatorname{SL}_2(q^3)\simeq
\operatorname{PGL}_2(q)$. We write $u$ as $u_1\cdot u_2$, where $u_1 \in
\operatorname{SL}_2(q^3)$, $u_2 \in \operatorname{SL}_2(q)$, and let $v_1, v_2$
be the images of
$u_1, u_2$ in $C/\operatorname{SL}_2(q)$ and
$C/\operatorname{SL}_2(q^3)$, respectively. Assume first that $q\equiv 1\pmod4$.
Then $\operatorname{PSL}_2(q)$ and
$\operatorname{PSL}_2(q^3)$ are strongly real. Therefore,  there exist
involutions 
$t_1 \in \operatorname{PSL}_2(q^3)$, $t_2 \in
\operatorname{PSL}_2(q)$ such that $v_1^{t_1}=v_1^{-1}$ and 
$v_2^{t_2}=v_2^{-1}$. Let $z_1$, $z_2$ belong to the preimages of
$t_1$, $t_2$ in $\operatorname{SL}_2(q^3)$ and
$\operatorname{SL}_2(q)$, respectively. Then $| z_1|=4=| z_2|$ and $z_1^2\in
Z(\operatorname{SL}_2(q^3))$, $z_2^2\in Z(\operatorname{SL}_2(q))$. It follows
that $z_1^2=z_2^2$ in $M$, whence $(z_1z_2)^2=e$. Thus $z_1z_2$ is an inverting
involution for $u$. 
Assume now that  $q\equiv 3\pmod4$. In this case there exist involutions $t_1\in
\operatorname{PGL}_2(q^3)\setminus \operatorname{PSL}_2(q^3)$, $t_2\in
\operatorname{PGL}_2(q)\setminus
\operatorname{PSL}_2(q)$ such that $v_1^{t_1}=v_1^{-1}$ and
$v_2^{t_2}=v_2^{-1}$. Furthermore, $t_1$ lies in a Cartan subgroup of
$\operatorname{PGL}_2(q^3)$, i.e., in a maximal torus of
$\operatorname{PGL}_2(q^3)$ of order $q^3-1$, while $t_2$ lies in a Cartan
subgroup of $\operatorname{PGL}_2(q)$, i.e., in a maximal torus of  
$\operatorname{PGL}_2(q)$ of order $q-1$. Let $T$ be a maximal torus of $C$ such
that its images under the natural homomorphisms 
$C\rightarrow C/\operatorname{SL}_2(q)$ and $C/\operatorname{SL}_2(q^3)$
contain elements $t_1$ and $t_2$, respectively.
Then  $| T|=(q^3-1)(q-1)$ and, by [17, Table~1.1], $T\simeq {\mathbb
Z}_{q^3-1}\times {\mathbb Z}_{q-1}$.
In particular,  $T$ does not contain elements of order $4$. Let $z$ be a
preimage of $t_1$ in $T$. We may assume that $z$ is a 
$2$-element; hence $z^2=e$. Moreover, since $t_1$ does not lie in
$\operatorname{PSL}_2(q^3)$,
we see that $z$ does not lie in $M$. Consider a natural homomorphism 
$\widetilde{\phantom{G}}:C\rightarrow C/\operatorname{SL}_2(q^3).$
Since $z\not\in M$, we obtain $\tilde{z}\not\in \operatorname{PSL}_2(q)$, and so
$t_2
\operatorname{PSL}_2(q)=\tilde{z} \operatorname{PSL}_2(q)$. Moreover,
$t_2,\tilde{z}\in
\widetilde{T}\simeq {\mathbb Z}_{q-1}$; hence $t_2=\tilde{z}$. Thus $z$ lies in
the 
preimage of $t_2$ as well. We obtain that $z$ is an
inverting involution for~$u$.

Let $s$ be such that its centralizer is conjugate either to the centralizer
of $s_5$ or to the centralizer of $s_{10}$. Then
$| C: (M\circ S)|=(2,q-1)$,  $M\simeq
\operatorname{SL}_2(q)$, and $S\simeq {\mathbb Z}_{q^3-\varepsilon}$, where
$\varepsilon=1$ if $C_G(s)$ is conjugate to $C_G(s_5)$ and $\varepsilon=-1$
if 
$C_G(s)$ is conjugate to $C_G(s_{10})$. Moreover, 
$C/S\simeq \operatorname{PGL}_2(q)$. We choose a maximal torus  $T$ of $C$ so
that $T\cap M$ is a Cartan subgroup of $M$. Since $M\simeq
\operatorname{SL}_2(q)$, we use matrices from $\operatorname{SL}_2(q)$ to write
elements of $M$, assuming that $T\cap M$ is the group of diagonal matrices.  By
Lemma~1, there exists an involution $x \in N(G,T)$ inverting each $t\in T$. In
particular,  $x$ inverts $s$ so $x$ normalizes $C_{\overline{G}}(s)$, i.e., 
$x\in N(G,C)$. Therefore, $x$ normalizes $\overline{S}$, and so it normalizes
$S$.
Put $C_0=\langle C, x\rangle$ and let $\widetilde{\phantom{C}}: C_0\rightarrow
C_0/S$ be the natural homomorphism. Since $M=O^{p'}(C)$ is characteristic in
$C$, $x$ induces an automorphism of $M$ of order  2. By [19, Lemma~2.3], 
$N(G,C)$ does not induce field automorphisms on $M$. Moreover, $x\not\in
\widehat{M}$, where $\widehat{M}$ is a group of inner-diagonal automorphisms of
${M}$, since $C/S\simeq
\operatorname{PGL}_2(q)\simeq \widehat{M}$. So
$C_0/S=\widetilde{C}_0\simeq \operatorname{PGL}_2(q)\times {\mathbb Z}_2.$ The
elements of $\widetilde{C}$ are written below as projective images of matrices
from 
$\operatorname{GL}_2(q)$. Up to conjugation in $\widetilde{C}_0$ we may assume
that 
$$\tilde{u}=\left[\begin{array}{rr} 1 & \alpha \\
0 & 1 \end{array}\right]$$
for some $\alpha \in {\mathbb F}_q$. Let
$\overline{T}\in\overline{G}$  be such that $T=\overline{T}_\sigma$.
Then $x$ normalizes $\overline{T}$. So, acting  by conjugation, 
$x$ leaves invariant the set of maximal unipotent subgroups of
$\overline{M}$, that are normalized by 
$\overline{T}\cap \overline{M}$. Furthermore, since $x$ is stable under 
$\sigma$; therefore,  $x$ normalizes the subgroups of $\sigma$-stable points of
these unipotent subgroups. Since
$\overline{M}=[\overline{C},\overline{C}]\simeq
\operatorname{SL}_2(\overline{{\mathbb F}}_q)$, there exists exactly two maximal
unipotent subgroups of $\overline{M}$ that are normalized by
$\overline{T}$:
one of them consists of upper-triangular matrices, another consists of
lower-triangular matrices. So 
$x$ either leaves these subgroups invariant or interchanges these subgroups.
Thus, either $\tilde{u}^{\tilde{x}}=\left[\begin{array}{rr} 1 &
\beta \\
0 & 1 \end{array}\right]$,  or
$\tilde{u}^{\tilde{x}}=\left[\begin{array}{rr} 1 & 0 \\
\beta & 1 \end{array}\right],$  for some $\beta \in {\mathbb F}_q$
Going back to elements $u$ and $x$ in $C_0$ and using the fact that $p$ is
coprime to
$|S|$ we derive that  $
u=\left(\begin{array}{rr} 1 & \alpha \\
0 & 1 \end{array}\right)$  and either
$u^x=\left(\begin{array}{rr} 1 & \beta \\
0 & 1 \end{array}\right),$ or $\left(\begin{array}{rr} 1 & 0 \\
\beta & 1 \end{array}\right).
$

Let $u^x=\left(\begin{array}{rr} 1 & \beta \\
0 & 1 \end{array}\right)$. Then there exists $\tilde{t}\in
\operatorname{PGL}_2(q)\cap
\widetilde{T}$ such that $${\left[\begin{array}{rr} 1 & \beta \\
0 & 1 \end{array}\right]}^{\tilde{t}}=\left[\begin{array}{rr} 1 & -\alpha \\
0 & 1 \end{array}\right].$$ Therefore,
$u^{xt}=\left(\begin{array}{rr} 1 & -\alpha \\
0 & 1 \end{array}\right)=u^{-1}$ and $(xt)^2=t^xt=t^{-1}t=e$.

Assume that $u^x=\left(\begin{array}{rr} 1 & 0 \\
\beta & 1 \end{array}\right)$. Then there exists $\tilde{t}\in
\operatorname{PGL}_2(q)\cap \widetilde{T}$ such that
$${\left[\begin{array}{rr} 1 & 0 \\
\beta & 1 \end{array}\right]}^{\tilde{t}}=\left[\begin{array}{rr} 1 & 0 \\
\alpha & 1 \end{array}\right].
$$
Therefore $u^{xt}=\left(\begin{array}{rr} 1 & 0 \\
\alpha & 1 \end{array}\right)=(u)^T$, where $^T$ denotes the transposition of a
matrix and $(xt)^2=t^xt=t^{-1}t=e$. Replacing $x$ by $xt$, we may assume
that $u^x=\left(\begin{array}{rr} 1 & 0 \\
\alpha & 1 \end{array}\right)=u^T$. Since  $| x|=2$, we also derive that
$(u^T)^x=u$. Set $z=\left(\begin{array}{rr} 0 & 1 \\
-1 & 0 \end{array}\right) \in \operatorname{SL}_2(q)\cap N(C_0, T)$. Then
$$u^{xz}=u^{-1}=u^{zx}.$$ Since
$| N(C_0, T)/T|=4$, $N(C_0, T)/T$ is abelian. Moreover, 
$x$ and $z$ lie in $N(C_0, T)$, and their images in $N(C_0, T)/T$ are
involutions;
hence  $N(C_0,T)/T\simeq
{\mathbb Z}_2\times{\mathbb Z}_2$ and $x$ normalizes $z (T\cap M)$, i.e.,
$z^{x}=zt$
for some $t\in T\cap M$. Therefore $xzt=zx$. As we noted above, both
$xz$ and $zx$ invert $u$, and so $t\in Z(M)$.
If $q$ is even then $Z(M)=\{e\}$, whence $x$
centralizes $\langle
z\rangle$. Therefore  $| xz|=2$, so $xz$ is an inverting involution. Assume
that  $q$ is odd. Then $|z|=4$,
$| Z(M)|=2$, and for $t\in Z(M)\setminus\{e\}$ the identity
$zt=z^{-1}$ holds. Thus either $z^x=z$, or
$z^x=z^{-1}$. We show that $z^x=z^{-1}$, whence
$(xz)^2=z^xz=z^{-1}z=e$ and $xz$ is an inverting involution. Let $Q$ be a
maximal torus of $C$ that contains $z$. Note that
$\tilde{x}$, $\widetilde{Q}$ lie in $C_{\widetilde{C}_0}(\tilde{z})$
and $C_{\widetilde{C}_0}(\tilde{z})\leq N(\widetilde{C}_0, \widetilde{Q})$.
Moreover, $\widetilde{C}_0=\operatorname{PGL}_2(q)\times \langle
\tilde{y}\rangle$ and
$\tilde{y} \in N({\widetilde{C}_0}, \widetilde{Q})$. Consider the cosets 
$\tilde{y}\widetilde{Q}$ and $\tilde{x}\widetilde{Q}$.
Suppose that these cosets coincide. Then 
$\tilde{x} \in \tilde{y}\widetilde{Q}$. Since $\widetilde{Q}$ is cyclic, it
contains a unique involution $\tilde{z}$. Therefore, 
$\tilde{y}\widetilde{Q}$ contains the two involutions
$\tilde{y}$, $\tilde{y}\tilde{z}$; hence $\tilde{x}=\tilde{y}$
or $\tilde{x}=\tilde{y}\tilde{z}$. The first equality is impossible, since 
$\tilde{y}$ centralizes $\operatorname{PGL}_2(q)$; and if 
$\tilde{x}=\tilde{y}\tilde{z}$, then
$\tilde{u}^{-1}=\tilde{u}^{\tilde{x}\tilde{z}}=\tilde{u}^{\tilde{y}\tilde{z}^2}
=\tilde{u}^{\tilde{z}^2}=\tilde{u}$, which is impossible.
Therefore $\tilde{y}\widetilde{Q}\neq\tilde{x}\widetilde{Q}$. By Lemma 1, there
exists an involution  $x' \in N(G, Q)$, inverting each element of
$Q$. We have $s \in Q$, and so
$x' \in C_0$ and $x, x' \in N(G,Q)\cap C_0$.
Since $\tilde{y}\widetilde{Q}\neq\tilde{x}\widetilde{Q}$,
$N(C_0,Q)/Q\simeq{\mathbb Z}_2\times{\mathbb Z}_2$ and $x,x'\not\in C$,
we see that
$\tilde{x}'\widetilde{Q}=\tilde{x}\widetilde{Q}$  and $xQ=x'Q$.
Therefore $x$ inverts each element of $Q$; in particular,~${z^x=z^{-1}}$.

Let  $s$ be such that its centralizer either is conjugate to the centralizer
of $s_3$, or is conjugate to the centralizer of  $s_{7}$. Then 
$| C: (M\circ S)|=(2,q-1)$,  $M\simeq
\operatorname{SL}_2(q^3)$,  and $S\simeq {\mathbb Z}_{q-\varepsilon}$, where
$\varepsilon=1$ if  $C_G(s)$ is conjugate to  $C_G(s_3)$ and
$\varepsilon=-1$ if
$C_G(s)$ is conjugate to $C_G(s_{7})$. Moreover, 
$C/S\simeq \operatorname{PGL}_2(q)$.
This case can settled by using exactly the same arguments as in the previous
case.

Assume that $C_G(s)$ is conjugate to $C_G(s_4)$. Then $M\simeq
\operatorname{SL}_3(q)$,
$S\simeq
{\mathbb Z}_{q^2+q+1}$. Moreover, if  $3$ divides $q-1$ then
$|C:M\circ S|=3$ and $C/S\simeq \operatorname{PGL}_3(q)$;
and if $3$ does not divide $q-1$ then  $C=
M\times S$ and $\operatorname{SL}_3(q)\simeq \operatorname{PGL}_3(q)$.
In both cases the proof is the same. Choose a maximal torus 
$T$ of $C$ so that $T\cap M$ is a Cartan subgroup of $M$. We identify elements
of $M$ with matrices of $\operatorname{SL}_3(q)$ and we assume that  $T\cap M$
is a subgroup of diagonal matrices under this identification. By Lemma~1,
there exists $x \in N(G,T)$ such that $x^2=e$ and $t^x=t^{-1}$ for every 
$t \in T$. Consider $C_0=\langle C,x\rangle$. By [19, Lemma~2.3], 
$N(G,C)$ does not induce field automorphisms on $M$. Since $x$ inverts each
element of a Cartan subgroup $T\cap M$ of $M$, we see that $x$ induces a
graph automorphism on $M$. Let
$\iota$ be a graph automorphism of  
$\operatorname{SL}_3(q)$ acting by  
$y\mapsto(y^{-1})^T$, where  ${}^T$ denotes the transposition of a
matrix. Then  $\iota$ normalizes $T\cap M$ and inverts each element from 
$T\cap M$. Hence, multiplying  $x$ by a suitable element of  $T\cap M$, we may
assume that  $x$ acts on $M$ in the same way as $\iota$. The element
$u$ is conjugate to its Jordan form in $C$, and so we may assume that 
$u=\left(\begin{array}{rrr} 1 & 1 & 0
\\ 0 & 1 & \alpha \\ 0 & 0 &
1\end{array}\right)$, where $\alpha \in
\{0, 1\}$. Let $u=\left(\begin{array}{rrr} 1 & 1 & 0 \\ 0 & 1 & 1 \\ 0 & 0
& 1\end{array}\right)$, and set
$z=\left(\begin{array}{rrr} 0 & 0 & -1 \\
0 & -1 & 0 \\ -1 & 0 &
0\end{array}\right) \in \operatorname{SL}_3(q)$. We have
$$
u^{xz}=((u^{-1})^T)^z
=\left(\left(\begin{array}{rrr} 1 & -1 & 1 \\ 0 & 1 & -1 \\
0 & 0 & 1\end{array}\right)^T\right)^z=\left(\begin{array}{rrr} 1 & 0 & 0 \\ -1
& 1 & 0 \\
1 & -1 & 1\end{array}\right)^z=\left(\begin{array}{rrr} 1 & -1 & 1 \\ 0 & 1 & -1
\\
0 & 0 &
1\end{array}\right)=u^{-1}.
$$ 
Therefore $xz$ is a sought involution, since 
$(xz)^2=z^xz=(z^{-1})^Tz=e$. Let
$u=\left(\begin{array}{rrr} 1 & 1 & 0 \\ 0 & 1 & 0 \\ 0 & 0 &
1\end{array}\right)$,  and set
$z=\left(\begin{array}{rrr} 0 & 1 & 0 \\
1 & 0 & 0 \\ 0 & 0 &
-1\end{array}\right) \in \operatorname{SL}_3(q)$. We have
$$u^{xz}=((u^{-1})^T)^z=\left(\left(\begin{array}{rrr} 1 & -1 & 0 \\ 0 & 1 & 0
\\
0 & 0 & 1\end{array}\right)^T\right)^z=\left(\begin{array}{rrr} 1 & 0 & 0 \\ -1
& 1 & 0 \\
0 & 0 & 1\end{array}\right)^z=\left(\begin{array}{rrr} 1 & -1 & 0 \\ 0 & 1 & 0
\\
0 & 0 & 1\end{array}\right)=u^{-1}.$$ Therefore $xz$ is a sought involution,
since  $(xz)^2=z^xz=(z^{-1})^Tz=e$.

Assume that  $C_G(s)$ is conjugate to $C_G(s_9)$. Then 
$M\simeq \operatorname{SU}_3(q)$, $S\simeq
{\mathbb Z}_{q^2-q+1}$. Moreover, if $3$ divides $(q+1)$ then 
$|C:M\circ S|=3$ and $C/S\simeq \operatorname{PGU}_3(q)$;
and  if  $3$ does not divide  $(q+1)$ then  $C=
M\times S$ and $\operatorname{SU}_3(q)\simeq \operatorname{PGU}_3(q)$.
In both cases the proof is the same. Choose a maximal torus 
$T$ of $C$ so that $T\cap M$ is a Cartan subgroup of $M$. By Lemma~1,
there exists $x \in N(G,T)$ such that $x^2=e$ and $t^x=t^{-1}$ for every 
$t \in T$. Again, by [19, Lemma~2.3], 
$N(G,C)$ does not induce field automorphisms on $M$. Let
$A=\left(\begin{array}{rrr} 0 & 0 & 1 \\ 0 & -1 & 0 \\ 1 & 0 &
0\end{array}\right)$ and let $\iota$ be an automorphism of
$\operatorname{SL}_3(q^2)$, acting by  $y\mapsto A(y^{-1})^TA$, where $^T$
denotes the transposition of a
matrix. Denote the automorphism of $\operatorname{SL}_3(q^2)$ that maps each
element of a matrix from  $\operatorname{SL}_3(q^2)$ into the power $q$ by  $f$.
In view of  [20, p.~268--270] we may assume that  $\operatorname{SU}_3(q)$
coincides with the set of  $\iota\circ f$-stable points.
We identify elements of $M$ with the set of $\iota\circ f$-stable points of
$\operatorname{SL}_3(q^2)$ and we assume that  $T\cap M$ is a group of diagonal
matrices under this identification. The restriction of $\iota$ on 
$\operatorname{SU}_3(q)$ we denote by the same symbol $\iota$. Then  
$\iota$ normalizes  $T\cap M$. So, multiplying $x$ by a suitable element of 
$T$, we may assume that $x$ acts on $M$ in the same way as $\iota$.
Up to conjugation in  $C$, each unipotent element $u$ of $M$ has the form
$u=\left(\begin{array}{rrr} 1 & \alpha & \beta
\\ 0 & 1 & \alpha^q \\ 0 & 0 &
1\end{array}\right)$, where $\alpha, \beta \in
F_{q^2}$ and $\beta + \beta^q=\alpha^{q+1}$. If $\alpha\neq 0$, there exists
an element $t \in T$ such that 
$u^t=\left(\begin{array}{rrr} 1 & 1 & \gamma \\ 0 & 1 & 1 \\
0 & 0 & 1\end{array}\right)$ for some $\gamma \in F_{q^2}$. So we may assume
that $u=\left(\begin{array}{rrr} 1 & 1 & \gamma
\\ 0 & 1 & 1 \\ 0 & 0 &
1\end{array}\right)$ and
$u^{-1}=\left(\begin{array}{rrr} 1 & -1 &
\gamma' \\ 0 & 1 & -1 \\ 0 & 0 &
1\end{array}\right)$ for some
$\gamma' \in F_{q^2}$. Put $z=\left(\begin{array}{rrr} -1 & 0 & 0 \\ 0 &
1 & 0 \\ 0 & 0 &
-1\end{array}\right) \in \operatorname{SU}_3(q)$. Then
$$u^{xz}=(A(u^{-1})^TA)^z=(A(u^{-1})^TA)^z=\left(\begin{array}{rrr} 1 & 1 &
\gamma' \\ 0 & 1 & 1 \\
0 & 0 & 1\end{array}\right)^z=\left(\begin{array}{rrr} 1 & -1 & \gamma' \\ 0 & 1
& -1 \\
0 & 0 & 1\end{array}\right)=u^{-1}.$$ For $\alpha=0$ the identity
$u^{xz}=u^{-1}$ is also true. Therefore, $xz$ is  a sought involution, since
$(xz)^2=z^xz=(z^{-1})^Tz=e$.

Theorem 1 and so Teorems 2 and 3 are proven.

\bigskip

\end{document}